\documentclass[reqno,12pt]{amsart}

\newtheorem{theorem}{Theorem}[section]
\newtheorem{lemma}[theorem]{Lemma}
\newtheorem{corollary}[theorem]{Corollary}

\theoremstyle{definition}
\newtheorem{assumption}[theorem]{Assumption}

\newtheorem{example}[theorem]{Example}
 
\theoremstyle{remark}
\newtheorem{remark}[theorem]{Remark}

\makeatletter
\def\dashint{\operatorname%
{\,\,\text{\bf--}\kern-.98em\DOTSI\intop\ilimits@\!\!}}
\makeatother

\newcommand\bR{\mathbb{R}}

\newcommand\cF{\mathcal{F}}

\newcommand\cL{\mathcal{L}}

\newcommand\cP{\mathcal{P}}
\newcommand\cS{\mathcal{S}}

\newcommand\frL{\mathfrak{L}}

\newcommand{\tr}{{\rm tr}\,}

\newcommand{\stref}[1]{\hbox{\rm\ref{#1}}}

 \newcommand{\mysection}[1]{\section{#1}
 \setcounter{equation}{0}}

\newcommand{\nlimsup}{\operatornamewithlimits{\overline{lim}}}

\begin{document}

\title[$L_{p}$-estimates]
{Some $L_{p}$-estimates for elliptic and parabolic
operators with measurable coefficients}
\author{N.V. Krylov}
\thanks{The work  was partially supported by
NSF Grant DMS-0653121}
\email{krylov@math.umn.edu}
\address{127 Vincent Hall, University of Minnesota,
 Minneapolis, MN, 55455}
 
 

\begin{abstract}
We consider linear elliptic and parabolic equations
with measurable coefficients
and prove two types of $L_{p}$-estimates
for their solutions,
 which were recently used
in the theory of fully nonlinear elliptic and parabolic
second order equations in \cite{DKL}. The first type
is an estimate of the $\gamma$th norm of the second-order derivatives,
where $\gamma\in(0,1)$, and the second type deals with
estimates of the resolvent operators in $L_{p}$ when
the first-order coefficients are summable to an appropriate power.
\end{abstract}

\maketitle

Let $d\geq1$ be an integer and let $\bR^{d}$ be a Euclidean space
of points $x=(x^{1},...,x^{d})$. 
Consider an operators $L$ of the form
\begin{equation}
                                                  \label{9.5.7}
L=\partial_{t}+a^{ij}(t,x)D_{ij}+b^{i}(t,x)D_{i}-c(t,x),
\end{equation}
where and below in the article the summation convention is enforced,
$$
D_{i}=\frac{\partial}{\partial x^{i}},\quad
D_{ij}=D_{i}D_{j},\quad \partial_{t}=\frac{\partial}{\partial t},
$$
  $a(t,x)=(a^{ij}(t,x))$ is a
uniformly nondegenerate and bounded matrix-valued,
$b(t,x)=(b^{i}(t,x))$ is an $\bR^{d}$-valued,
and $c(t,x)$ is a real-valued measurable functions
defined on $\bR^{d+1}=\{(t,x):t\in\bR,x\in\bR^{d}\}$.

In this article we are going to discuss two types of estimates
for operators like $L$, which were recently used
in the theory of fully nonlinear elliptic and parabolic
second order equations in \cite{DKL}.

 The first type
(see Theorems \ref{theorem 9.3.2} and \ref{theorem 8.16.1})
is about the possibility to estimate the integrals of $|D^{2}u|^{\gamma}$
with some $\gamma\in(0,1)$ through the $L_{p}$-norm
of $Lu$ and the sup norm of $u$, where $D^{2}u$
is the Hessian of $u$. This seemingly very weak estimate,
discovered for elliptic equations by F.H. Lin, recently played
a crucial role in the theory of fully nonlinear elliptic
and parabolic equations with VMO ``coefficients'' (see
\cite{DKL}). In \cite{DKL} we use a result stated in \cite{Kr10}
without proof. Even though the proof is not difficult
it is still worth presenting it with all details especially
because on our way we obtain some new nontrivial information
such as Lemma \ref{lemma 8.16.1} or its probabilistic counterpart
Theorem \ref{theorem 9.4.1}. One more point worth mentioning is
that unlike F.H. Lin, who used a rather delicate reversed H\"older's inequality
which was proved an the basis of Gehring's lemma, we are using a
basic result  of Krylov-Safonov, which provided the foundation
of the theory of fully nonlinear elliptic and parabolic second-order
equations. In fact, we need its version
obtained in \cite{Kr10} just by analyzing the corresponding arguments
in \cite{Kr87}.
To obtain the above mentioned
estimate we assume that $b$ and $c$ are bounded. Similar estimates
we give for $|Du|$, where $Du$ is the gradient of $u$.

The second type of results
deals with estimates of the $L_{p}$-norm of $ \mu u$ through the
$L_{p}$-norm of $\mu u-Lu$ if $\mu>0$
with a constant independent of $\mu$
if $\mu$ is large (see Theorems \ref{theorem
8.25.1} and \ref{theorem 9.10.1}). As we have noted
 these theorems are also
used in \cite{DKL} in particular cases when the drift coefficients
are bounded. However, even in this case we could not find  a direct
reference to the result we needed and, therefore, our explanation
in \cite{DKL}
contains the words such as ``by analyzing the proof...''.
Here we prove the corresponding result with all details 
and also give its generalization for the case in which
$b$ is in $L_{q}$ with an appropriate $q\leq p$.

\mysection{Estimates of $|D^{2}u|$}

Fix a $\delta\in(0,1)$ and introduce
$\cS_{\delta}$ as the set of symmetric $d\times d$-matrices
$a=(a^{ij})$ such that for any $\xi\in\bR^{d}$ we have
$$
\delta|\xi|^{2}\leq a^{ij}\xi^{i}\xi^{j}\leq\delta^{-1}|\xi|^{2}.
$$
For constants $K\geq0$
denote by $\frL_{\delta,K}$ the set of operators $L$ of type
\eqref{9.5.7},
when $a(t,x)=(a^{ij}(t,x))$ is  $\cS_{\delta}$-valued and
$b$ and $c$ are such that
$$
|b |+c \leq K,\quad c \geq0.
$$
Let $\frL^{0}_{\delta,K}$ be a subset of $\frL_{\delta,K}$ consisting of
operators with infinitely differentiable coefficients.

For $\rho,r>0$ introduce 
$$
B_{r}=\{x\in\bR^{d}:|x|<r\},
$$
$$
C_{\rho,r}=(0,\rho)\times B_{r},\quad
\partial'C_{\rho,r}=\big([0,\rho]\times\partial B_{r}\big)
\cup\big(\{\rho\}\times B_{r}\big),\quad C_{r}=C_{r^{2},r},
$$
$$
C_{\rho,r}(t,x)=(t,x)+C_{\rho,r},\quad
\partial'C_{\rho,r}(t,x)=(t,x)+\partial'C_{\rho,r}.
$$

Our first goal in this section is to prove the following
parabolic version of the main result of
\cite{FHL} by  F.H. Lin.

\begin{theorem}
                                             \label{theorem 8.11.1}

There are   constants $\gamma\in(0,1]$ and $N$,
depending only on $\delta$, $K$, and $d$, such that
for any $L
\in\frL_{\delta,K}$ and $u\in  W^{1,2}_{d+1,loc}(C_{2,1})\cap C(\bar{C}_{2,1})  $
 we have
\begin{equation}
                                                    \label{8.11.1}
 \int_{C_{1,1} (1,0) }|D^{2}u|^{\gamma}
\,dxdt \leq N\big(\int_{C_{2,1}}|Lu|^{d+1}\,dxdt\big)^{\gamma/(d+1)}
 +N\sup_{\partial'C_{2,1}}|u|^{\gamma} .
\end{equation}
\end{theorem}
This theorem is stated as Corollary 4.2 in \cite{Kr10} but no 
proof is given there. We fill this gap in this article.

The following theorem is proved in \cite{Kr10}.

\begin{theorem}
                                           \label{theorem 12.21.1}
Let $u\in W^{1,2}_{d+1}(C_{2,1})$ and assume that
$u\geq0$ on $\partial'C_{2,1}$ and there exists
an operator $L\in\frL_{\delta,K}$ such that $Lu\leq0$ in $C_{2,1}$.
Then there exist  constants $\gamma=\gamma(\delta,d,K)\in(0,1]$
and $N =N (\delta,d,K)$
such that for any $\lambda>0$
\begin{equation}
                                                    \label{12.21.2}
|C_{1,1}(1,0)\cap\{(t,x):-Lu(t,x)\geq\lambda\}| 
\leq N \lambda^{-\gamma}u^{\gamma}(0,0).
\end{equation}

\end{theorem}

\begin{corollary}
                                               \label{corollary 8.11.1}
Under the conditions of Theorem \stref{theorem 12.21.1}
for any $\gamma'\in(0,\gamma)$
we have
\begin{equation}
                                                    \label{8.11.2}
 \int_{C_{1,1} (1,0) }|Lu|^{\gamma'}
\,dxdt \leq Nu^{\gamma'}(0,0) ,
\end{equation}
where $N =N (\delta,d,K,\gamma')$.

\end{corollary}

Indeed, the left-hand side of \eqref{8.11.2} equals
$$
\int_{0}^{\infty}
|C_{1,1}(1,0)\cap\{-Lu\geq\lambda^{1/\gamma'}\}|\,d\lambda
$$
$$
\leq\int_{0}^{\mu}|C_{1,1}(1,0) \}|\,d\lambda
+Nu^{\gamma}(0,0)\int_{\mu}^{\infty}\lambda^{-\gamma/\gamma'}\,d\lambda,
$$
where $\mu=u^{\gamma'}(0,0)$. Upon computing the last integral we arrive at
\eqref{8.11.2}.

For elliptic operators we have the following version of
Theorem \stref{theorem 12.21.1}.
\begin{theorem}
                                               \label{theorem 9.3.1}
Let $u\in W^{2}_{d}(B_{1})$ and assume that
$u\geq0$ on $\partial B_{1}$ and there exists
an operator $L\in\frL_{\delta,K}$ with coefficients independent of $t$
such that
$Lu\leq0$ in
$B_{1}$. Then there exist  constants $\gamma=\gamma(\delta,d,K)\in(0,1]$
and $N =N (\delta,d,K)$
such that for any $\lambda>0$
\begin{equation}
                                                    \label{9.3.1}
|B_{1} \cap\{x:-Lu(x)\geq\lambda\}| 
\leq N \lambda^{-\gamma}u^{\gamma}(0).
\end{equation}

\end{theorem}

Proof. First assume that $u\in C^{2}(\bar{B}_{1})$ and define a function
$v=v(t,x)$ by $v(t,x)=u(x)$. By the maximum principle 
$u\geq0$ in $B_{1}$. Therefore $v$ satisfies the assumptions
of Theorem \stref{theorem 12.21.1} and \eqref{9.3.1}
in this particular case follows from
\eqref{12.21.2}.

In the general case, introduce $f=-Lu$ and find a sequence of
operators $L_{n}\in\frL_{\delta,K}$, $n=1,2,...$, with smooth coefficients
converging (a.e.) to the corresponding coefficients of $L$. Also
let $f_{n}\in C^{1}(\bar{B}_{1})$, $n=1,2,...$, be a sequence
of nonpositive functions such that $f_{n}\to f$ in $L_{d}(B_{1})$.
Define $u_{n}\in C^{2}(\bar{B}_{1})$
as unique solutions of equations $L_{n}u_{n}=-f_{n}$
with zero boundary condition. Since, $u_{n}\leq u$ on $\partial B_{1}$
and $L_{n}(u_{n}-u)=-f_{n}+f+(L-L_{n})u\to0$ in $L_{d}(D_{1})$,
we have that 
$$
\nlimsup_{n\to\infty}u_{n}\leq u
$$
in $\bar{B}_{1}$ owing to the Alexandrov estimate. Now we recall that
the convergence almost everywhere implies the 
convergence in distribution and conclude that
$$
F(\lambda):=|B_{1} \cap\{x:f(x)\geq\lambda\}|=
\lim_{n\to\infty}|B_{1} \cap\{x:f_{n}(x)\geq\lambda\}|
$$
$$
\leq N\lambda^{-\gamma}\nlimsup_{n\to\infty} u_{n}^{\gamma}(0)
\leq N\lambda^{-\gamma}u^{\gamma}(0)
$$
at all $\lambda>0$ at which $F(\lambda)$ is continuous. Since the
right-hand side of \eqref{9.3.1} is continuous in $\lambda$
and the left-hand side is right continuous, we have 
\eqref{9.3.1} for all $\lambda>0$ and the theorem is proved.

As in the case of Corollary \ref{corollary 8.11.1}
we have the following.

\begin{corollary}
                                               \label{corollary 4.25.1}
Under the conditions of Theorem \stref{theorem 9.3.1}
for any $\gamma'\in(0,\gamma)$
we have
$$
 \int_{B_{1} }|Lu|^{\gamma'}
\,dx  \leq Nu^{\gamma'}( 0) ,
$$
where $N =N (\delta,d,K,\gamma')$.

\end{corollary}

Here is a useful generalization of Corollary \stref{corollary 8.11.1}.
\begin{lemma}
                                                   \label{lemma 8.16.1}

There are   constants $\gamma\in(0,1]$ and $N$,
depending only on $\delta$, $K$, and $d$, such that,
if $L
\in\frL_{\delta,K}$ and $u\in  W^{1,2}_{d+1,loc}(C_{2,1})\cap C(\bar{C}_{2,1})  $
and $Lu=g-f$ in $C_{2,1}$ with $f\geq0$, then
 we have
$$
 \int_{C_{1,1} (1,0) }|f|^{\gamma}
\,dxdt \leq N\big(\int_{C_{2,1}}|g_{+}|^{d+1}\,dxdt\big)^{\gamma/(d+1)}
$$
\begin{equation}
                                                    \label{8.16.4}
 +N |u_{+}(0,0)|^{\gamma} +N\sup_{\partial'C_{2,1}}|u_{-}|^{\gamma}.
\end{equation}
\end{lemma}
Proof. 
First we reduce the general case to the one in which
$u\in W^{1,2}_{d+1}(C_{2,1})$.
To do that we introduce ``shifted and dilated'' $u$,
that is for $\varepsilon\in[0, 1)$ we define
$$
u_{\varepsilon}(t,x)=u( \varepsilon^{2} t+ 1-\varepsilon^{2},
 \varepsilon x) ).
$$
Obviously,
$u_{\varepsilon}\in   C(\bar C_{1})\cap  W^{1,2}_{d+1 }(C_{2,1})$.
We also modify the coefficients of $L$ in such a way that
the new $g$ and $f$ are just shifted and dilated original $g$
and $f$, respectively,
times $\varepsilon^{2}$.
If \eqref{8.16.4} holds for $u_{\varepsilon}$,
then we obtain it as is by letting $\varepsilon\uparrow1$
by the monotone convergence theorem owing to the continuity
of $u$ in $\bar{C}_{2,1}$. By the way, we do not assume that 
the integral in the right-hand side of \eqref{8.16.4} is finite.
Thus indeed we may concentrate
on $u\in W^{1,2}_{d+1}(C_{2,1})$.

Then observe that if
the integral in the right-hand side of \eqref{8.16.4} is infinite,
we have nothing to prove. Therefore, we may assume that it is finite.
Then
$g\in L_{d+1}(C_{2,1})$ since $g_{+}\geq g\geq Lu$. It follows that
$f\in L_{d+1}(C_{2,1})$ as well.

Now take an operator $L'\in\cL^{0}_{\delta}$ and introduce
$$
f'=f+[(L'-L)u]_{-},\quad g'=g+[(L'-L)u]_{+},
$$
 so that
$L'u=g'-f'$ and $f',g'\in L_{d+1}(C_{2,1})$ and $f'\geq0$.
Obviously, if \eqref{8.16.4} were true with $f',g'$ in place of
$f,g$ for any $L'$, then by approximating $L$ by operators $L'$
we would obtain \eqref{8.16.4} in its original form.

Therefore, in the rest of the proof without losing generality
we assume that $L \in\cL^{0}_{\delta}$ and introduce functions $v$
and $w$ as $W^{1,2}_{d+1}(C_{2,1})$ solutions of
$$
Lv=-f,\quad Lw=-g
$$
with zero condition on $\partial'C_{2,1}$ for $v$ and with
condition $w=-u$ on  $\partial'C_{2,1}$. The existence and uniqueness
of $v$ and $w$ is a classical result (see, for instance, \cite{Li}).

Clearly, $v=u+w$ and by the maximum principle $v\geq0$. By 
Corollary \stref{corollary 8.11.1}, for an appropriate $\gamma$,
 the left-hand
side of
\eqref{8.16.4} is less than a constant times
$$
v^{\gamma}(0,0)\leq u_{+}^{\gamma}(0,0)+w_{+}^{\gamma}(0,0).
$$
After that it only remains to use the parabolic Alexanrdov
estimate. The lemma is proved.

For elliptic operators Lemma \stref{lemma 8.16.1} becomes the following.
\begin{lemma}
                                                   \label{lemma 9.3.2}

There are   constants $\gamma\in(0,1]$ and $N$,
depending only on $\delta$, $K$, and $d$, such that,
if $L
\in\frL_{\delta,K}$ has coefficients
independent of $t$ and $u\in 
W^{2}_{d,loc}(B_{1})\cap C(\bar{B}_{1})  $ and $Lu=g-f$ in
$B_{1}$ with $f\geq0$, then
 we have
$$
 \int_{B_{1}}|f|^{\gamma}
\,dx \leq N\big(\int_{B_{1}}|g_{+}|^{d }\,dx \big)^{\gamma/d}
 +N |u_{+}( 0)|^{\gamma} +N\sup_{\partial B_{1}}|u_{-}|^{\gamma}.
$$
\end{lemma}

The proof is based on Corollary \ref{corollary 4.25.1}
and consists of  repeating
 the proof of Lemma \ref{lemma 8.16.1} with obvious changes.
Of course, at the last step 
one applies the original Alexandrov estimate rather than
its parabolic version.

{\bf Proof of Theorem \ref{theorem 8.11.1}}.
Introduce $h=Lu$, take an operator $L'\in\frL_{\delta/2,K}$,
and observe that
$$
L'u=g-f,\quad g:=h+2[(L'-L)u]_{+},\quad f:=|(L'-L)u|.
$$
According to \eqref{8.16.4} and the parabolic Alexanrdov
estimate  
$$
\int_{C_{1,1} (1,0) }|(L'-L)u|^{\gamma }
\,dxdt \leq  N\|[(L'-L)u]_{+}\|_{L_{d+1}(C_{2,1})}^{\gamma }
$$
\begin{equation}
                                                    \label{8.11.3}
+N\|h\|_{L_{d+1}(C_{2,1})}^{\gamma }
+N\sup_{\partial'C_{2,1}}|u|^{\gamma }.
\end{equation}

We now use the arbitrariness of $L'$. Obviously, there exists
an $\varepsilon=\varepsilon(\delta,d)>0$ and an operator $L'\in\cL_{\delta/2}$
with lower order coefficients coinciding with the ones of $L$ and such that
$$
L'u=Lu-\varepsilon|D^{2}u| .
$$
With such an operator \eqref{8.11.3} becomes \eqref{8.11.1}.
The theorem is proved.

The reader understands that the following result
is obtained by mimicking the proof
of Theorem \stref{theorem 8.11.1} and using Lemma
\stref{lemma 9.3.2} instead of Lemma~\stref{lemma 8.16.1}.

\begin{theorem}
                                             \label{theorem 9.3.2} 

There are   constants $\gamma\in(0,1]$ and $N$,
depending only on $\delta$, $K$, and $d$, such that 
for any $L
\in\frL_{\delta,K}$ with the coefficients
independent of $t$ and $u\in 
W^{2}_{d,loc}(B_{1})\cap C(\bar B_{1})  $
 we have
\begin{equation}
                                       \label{9.5.6}
 \int_{B_{1}}|D^{2}u|^{\gamma}
\,dx \leq N\big(\int_{B_{1}}|Lu|^{d}\,dxdt\big)^{\gamma/d}
 +N\sup_{\partial B_{1}}|u|^{\gamma} .
\end{equation}
\end{theorem}

Next result is 
stronger than Theorem \stref{theorem 8.11.1} and looks like the right
parabolic counterpart of
Theorem \stref{theorem 9.3.2}.
It is proved in \cite{DKL} on the basis of Theorem \stref{theorem 8.11.1}.
We give it with a proof just for completeness.

\begin{theorem}
                       \label{theorem 8.16.1}
Let  
$u\in   C(\bar C_{1})\cap  W^{1,2}_{d+1,loc}(C_{1})$. Then there are
constants $\gamma\in(0, 1]$
 and $N$,
depending only on
$\delta, d$, and $K$, such that
for any $L\in  \frL_{\delta,K} $ we have
\begin{equation}
                                                          \label{8.11.01}
\int_{C_{1}} |D^2u|^\gamma \, dx \,  dt
 \leq  N \sup_{\partial'C_{1}}
 |u|^\gamma + N
\left(\int_{C_{1}}|  Lu|^{d+1} \,
dx\, dt\right)^{\gamma/{(d+1)}}.
\end{equation}  
\end{theorem}

Proof.   
First as in the proof of Lemma \ref{lemma 8.16.1}
one reduces the general situation  to the one in which
$u\in W^{1,2}_{d+1}(C_{1})$.

  Then
we may also assume that the coefficients  of $L$
are infinitely differentiable in $\bR^{d+1}$. Now
set $f=  Lu$ in $C_{1}$ and extend $f(t,x)$
for $t\leq 0$ as zero. Also set $u(t,x)=u(-t,x)$
for $t\leq0$. Observe that the new $u$ belongs to
$W^{1,2}_{d+1 }((-1,1)\times B_{1})$.
After that
define $v(t,x)$ as a unique
 $W^{1,2}_{d+1,\text{loc}}((-1,1)\times B_{1})\cap C([-1,1]\times \bar B_1)$
solution of
 $ Lv=f$ with  terminal
and lateral conditions being $u$. The existence and uniqueness
of such a solution is a classical result
(see, for instance,   Theorem 7.17 of \cite{Li}).
By uniqueness $v=u$ in $C_{1}$, so that
owing to Theorem \stref{theorem 8.11.1}, 
$$
\int_{ C_1 } |D^2u|^\gamma \, dx\, dt
=\int_{ C_1 } |D^2 v|^\gamma \, dx\, dt\leq N
\left(\int_{(-1,1)\times B_{1}}|f|^{d+1}\,dx\,dt
\right)^{\gamma/(d+1)}
$$
$$
+N\sup_{\partial'(-1,1)\times B_{1}} |v|^\gamma 
=N
\left(\int_{C_{1}}|f|^{d+1}\,dx\,dt
\right)^{\gamma/(d+1)}
+N\sup_{\partial'C_{1}} |u|^\gamma .
$$
The theorem is proved.

\mysection{Estimating $|Du|$}
 
\begin{lemma}
                                                   \label{lemma 8.16.04}

There are   constants $\gamma\in(0,1]$ and $N$,
depending only on $\delta$, $K$, and $d$, such that,
if $L
\in\frL_{\delta,K}$ and $u\in  W^{1,2}_{d+1,loc}(C_{2,1})\cap C(\bar{C}_{2,1})$,
then
 we have
\begin{equation}
                                                    \label{8.16.5}
\int_{C_{1,1} (1,0) }|Du|^{\gamma}
\,dxdt \leq
N\big(\int_{C_{2,1}}|Lu|^{d+1}\,dxdt\big)^{\gamma/(d+1)}
+N\sup_{\partial'C_{2,1}}|u|^{\gamma}.
\end{equation}
\end{lemma}

Proof. It certainly suffices to concentrate on smooth $u$.
In that case
observe that
\begin{equation}
                                                            \label{9.3.5}
L(-u^{2})=g-f,\quad g:=-2uLu-cu^{2},\quad f=2a^{ij}(D_{i}u)D_{j}u.
\end{equation}

By Lemma \stref{lemma 8.16.1} with an appropriate $\gamma$    
$$
\int_{C_{2,1}}|Du|^{2\gamma}\,dx\,dt\leq N 
\sup_{C_{2,1}}|u|^{\gamma} \big(\int_{C_{2,1}}|Lu|^{d+1}\,dx\,dt
\big)^{\gamma/(d+1)}
+N\sup_{C_{2,1}}|u|^{2\gamma}.
$$
After that it only remains to use 
Jensen's inequality and again the parabolic Alexandrov estimate.
The lemma is proved.

We also have \eqref{9.3.5} for elliptic operators.
Therefore, as above, Lemma \stref{lemma 9.3.2} yields
\begin{theorem}
                                                   \label{theorem 9.3.5}

There are   constants $\gamma\in(0,1]$ and $N$,
depending only on $\delta$, $K$, and $d$, such that,
if $L
\in\frL_{\delta,K}$ has the coefficients
independent of $t$ and $u\in  W^{2}_{d
,loc}(B_{1})\cap  C(\bar{B}_{1})$,
then
 we have
$$
 \int_{B_{1}}|Du|^{\gamma}
\,dx  \leq
N\sup_{\partial B_{1}}|u|^{\gamma}
 +N\big(\int_{B_{1}}|Lu|^{d }\,dx
\big)^{\gamma/d}.
$$
\end{theorem}

Here is our estimate of $Du$ in the parabolic case.

\begin{theorem}
Let  
$u\in   C(\bar C_{1})\cap  W^{1,2}_{d+1,loc}(C_{1})$.
 Then there are constants $\gamma\in(0, 1]$ and $N$,
depending only on
$\delta$, $K$, and $d$, such that
for any $L\in  \cL_{\delta,K} $ we have
$$
\int_{C_{1}} |Du|^{ \gamma} \, dx \,  dt
 \leq   N 
\left(\int_{C_{1}}|  Lu|^{d+1} \,
dx\, dt\right)^{ \gamma/{(d+1)}}+
N \sup_{\partial'C_{1}}
 |u|^{ \gamma} .
$$
\end{theorem}

This theorem is derived from Lemma \stref{lemma 8.16.04}
in the same way as Theorem \stref{theorem 8.16.1}
is derived from Theorem \stref{theorem 8.11.1}.

\mysection{Probabilistic versions}

Let $(\Omega,\cF,P)$  be
 a complete probability space
endowed with an increasing filtration of
$\sigma$-fields $\cF_{t}\subset\cF$, $t\geq0$, each of which is complete
with respect to $P,\cF$. By $\cP$ we denote
the predictable $\sigma$-field on $\Omega\times(0,\infty)$
generated by $\{\cF_{t},t\geq0\}$. Let $w_{t}$,
$t\geq0$, be
 a $d_{1}$-dimensional $\cF_{t}$-Wiener process   on
on $\Omega$, where $d_{1}\geq d$ is an integer.  Assume that  on
$\Omega\times(0,\infty)$ we are given $\cP$-measurable functions
$\sigma_{t}=\sigma_{t}(\omega)$ and $b_{t}=b_{t}(\omega)$ with values in
the set of $d\times d_{1}$-matrices and $\bR^{d}$, respectively.
Suppose that $a_{t}:=(1/2)\sigma_{t}\sigma^{*}_{t}\in\cS_{\delta}$
and $|b_{t}|\leq K$ for all $(\omega,t)$, where $K$ 
and $\delta$ are fixed constants.

 \begin{theorem}
                                                 \label{theorem 9.4.1}
Introduce
$$
x_{t}=\int_{0}^{t}\sigma_{s}\,dw_{s}
+\int_{0}^{t}b_{s}\,ds,
$$
$$
\tau=\inf\{t\geq0:(t,x_{t})\not\in C_{2,1}\}.
$$
Let $f(t,x)$, $(t,x)\in(0,\infty)\times\bR^{d}$, be a nonnegative Borel
function such that
$f(t,x)=0$ for $t\leq 1$. Then
\begin{equation}
                                                      \label{9.4.1}
\big(\int_{C_{2,1}}f^{\gamma}\,dxdt\big)^{1/\gamma}
\leq NE\int_{0}^{\tau}f(t,x_{t})\,dt,
\end{equation}
where $\gamma=\gamma(\delta,d,K)>0$ and $N=N(\delta,d,K)<\infty$.

\end{theorem}

Proof. First assume that $f$ is infinitely differentiable
in $(t,x)$.
Consider the following Bellman's equation:
\begin{equation}
                                                    \label{9.5.1}
\partial_{t}u+\inf_{a\in\cS_{\delta},|b|\leq K}
\big[a^{ij}D_{ij}u+b^{i}D_{i}u+f\big]=0
\end{equation}
in $C_{2,1}$ with zero boundary data on $\partial'C_{2,1}$.
By Theorem 6.4.1 of \cite{Kr87} this problem has a unique solution
bounded and continuous in $\bar{C}_{2,1}$ and having
bounded and continuous in $C_{2,1}$ derivatives $\partial_{t}u$,
$Du$, and $D^{2}u$.
Actually, to apply Theorem 6.4.1 of \cite{Kr87} directly
we need to have the term $-u$ in the left-hand side of \eqref{9.5.1}.
However, this is easily achieved by introducing a new function
$v$ such that $u=e^{-t}v$.
 By the maximum principle $u\geq0$.

Obviously,
\begin{equation}
                                                          \label{9.5.2}
\partial_{t}u(t,x_{t})+a^{ij}_{t}D_{ij}u(t,x_{t})
+b^{i}_{t}D_{i}u(t,x_{t})+f(t,x_{t})\geq 0,
\end{equation}
for $t<\tau$. Furthermore, it is easy to see that  there exists an
operator
$L\in\frL_{\delta,K}$ such that $Lu=-f$, so that by Lemma 
\ref{lemma 8.16.1}
\begin{equation}
                                                           \label{9.5.3}
 \int_{C_{1,1} (1,0) }|f|^{\gamma}
\,dxdt\leq Nu^{\gamma}(0,0).
\end{equation}

 Due to \eqref{9.5.2}, by It\^o's formula
$$
u(t\wedge\tau,x_{t\wedge \tau})=
u(0,0)+m_{t}
$$
$$
+\int_{0}^{t\wedge\tau}\big[
\partial_{t}u(t,x_{t})+a^{ij}_{t}D_{ij}u(t,x_{t})
+b^{i}_{t}D_{i}u(t,x_{t})\big]\,dt
$$
$$
\geq u(0,0)-\int_{0}^{t\wedge\tau}f(s,x_{s})\,ds+m_{t},
$$
where $m_{t}$ is a martingale. Upon plugging in $t=4$,
observing that $4\wedge\tau=\tau$ and $u(\tau,x_{\tau})=0$,
and taking the expectations of the extreme terms we obtain
$$
E\int_{0}^{\tau}f(t,x_{t})\,dt\geq u(0,0).
$$
After that to prove \eqref{9.4.1}
for infinitely differentiable $f$, it only remains to use \eqref{9.5.3}.

The parabolic Alexandrov estimate in probabilistic terms (see, for
instance, Theorem 2 of
\cite{Kr74} or Theorem 2.2.4 of \cite{Kr77}) implies that
for any Borel $g\geq 0$
$$
E\int_{0}^{\tau}g(t,x_{t})\,dt\leq N\|g\|_{L_{d+1}(C_{2,1})},
$$
where $N=N(d,\delta,K)$. This easily allows us to extend
\eqref{9.4.1} to the set of bounded Borel $f$
(vanishing for $t\leq1$). Finally, applying
the monotone convergence theorem
we get the desired result.
The theorem is proved.

We have derived Theorem \ref{theorem 9.4.1}
from Lemma \ref{lemma 8.16.1}, but actually 
Theorem \ref{theorem 9.4.1} is equivalent to
Lemma \ref{lemma 8.16.1}

In probabilistic terms Lemma \ref{lemma 9.3.2} means the following.
\begin{theorem}
                                      \label{theorem 4.29.1}
There exist constants $\gamma\in(0,1)$ and $N\in(0,\infty)$
depending only on $\delta,K$, and $d$, such that
if $f(x)$ is a nonnegative function on $B_{1}$ and
$$
\tau=\inf\{t\geq0:|x_{t}|=1\},
$$
then
\begin{equation}
                                                   \label{4.29.5}
\big(\int_{B_{1}}f^{\gamma}\,dx\big)^{1/\gamma}
\leq NE\int_{0}^{\tau}e^{-Kt}f(x_{t})\,dt
\leq NE\int_{0}^{\tau}f(x_{t})\,dt.
\end{equation}
\end{theorem}

We leave it to the reader to follow closely the proof
of Theorem \ref{theorem 9.4.1} by using the corresponding
results for elliptic equations from \cite{Kr87}.

The probabilistic versions of Lemmas \ref{lemma 8.16.1}
and \ref{lemma 9.3.2} allow  one to 
give   different proofs of Theorems \ref{theorem 9.3.2}
and \ref{theorem 8.16.1}. We only show this on the example
of Theorem \ref{theorem 9.3.2}.

{\bf Proof of Theorem \ref{theorem 9.3.2}}. First as in the proof of
Lemma \ref{lemma 8.16.1} we may assume that $u\in W^{2}_{d}(B_{1})$.
Then
we find an $\varepsilon=\varepsilon(\delta,d)>0$ and an
operator
$L'\in\frL_{\delta/2,K}$ with coefficients independent of $t$
such that
$$
L'u=Lu-\varepsilon |D^{2}u|.
$$

Let $a=a(x)=(a^{ij}(x))$, $b=b(x)=(b^{i}(x))$, and $c=c(x)$
be the coefficients of $L'$. Define $\sigma=\sqrt{2a}$.
One knows (see, for instance, \cite{Kr69} or
\cite{Kr77}) that there always exist  
$(\Omega,\cF,P)$, $\cF_{t}$, and $w_{t}$ as in the beginning
of the section,  
and there exists an $\cF_{t}$-adapted
continuous $\bR^{d}$-valued process $x_{t}$, $t\geq0$,
on $\Omega$ such that with probability one for all $t\geq0$
$$
x_{t}=\int_{0}^{t}\sigma(x_{s})\,dw_{s}
+\int_{0}^{t}b(x_{s})\,ds.
$$
From \cite{Kr74} (see the comments
after Theorem 3 there and see Theorem 4 of
\cite{Kr69}) or \cite{Kr77} we know that It\^o's formula
is applicable to 
$$
u(x_{t\wedge\tau})\exp\big[-\int_{0}^{t\wedge\tau}c(x_{s})\,ds\big].
$$
Therefore,
$$
u(0)=Eu(x_{\tau})\exp\big[-\int_{0}^{\tau}c(x_{s})\,ds\big]
$$
$$ 
-E\int_{0}^{\tau}L'u(x_{t})
\exp\big[-\int_{0}^{t }c(x_{s})\,ds\big] \,dt.
$$
By using the probabilistic version of the Alexandrov estimate
and the fact that $0\leq c\leq K$
we conclude
$$
\varepsilon  
E\int_{0}^{\tau}|D^{2}u|(x_{t})e^{-Kt}\,dt\leq
\varepsilon
E\int_{0}^{\tau}|D^{2}u|(x_{t})
\exp\big[-\int_{0}^{t }c(x_{s})\,ds\big] \,dt 
$$
$$
=u(0)+E\int_{0}^{\tau}Lu(x_{t})
\exp\big[-\int_{0}^{t }c(x_{s})\,ds\big] \,dt
$$
$$
-Eu(x_{\tau})\exp\big[-\int_{0}^{\tau}c(x_{s})\,ds\big] \leq
u(0)+N\|Lu\|_{L_{d}(B_{1})}+\sup_{\partial B_{1}}|u|,
$$
so that by Theorem \ref{theorem 4.29.1}
\begin{equation}
                                               \label{9.5.5}
\big(\int_{B_{1}}|D^{2}u|^{\gamma}\,dx\big)^{1/\gamma}
\leq N|u(0)|+N\|Lu\|_{L_{d}(B_{1})}+N\sup_{\partial B_{1}}|u|.
\end{equation}

Now observe that the above argument is applicable
for $\varepsilon=0$ and $L'=L$ in which case we get
$$
0=u(0)+E\int_{0}^{\tau}Lu(x_{t})
\exp\big[-\int_{0}^{t }c(x_{s})\,ds\big] \,dt
$$
$$
-Eu(x_{\tau})\exp\big[-\int_{0}^{\tau}c(x_{s})\,ds\big],
$$
$$
|u(0)|\leq
\big|E\int_{0}^{\tau}Lu(x_{t})
\exp\big[-\int_{0}^{t }c(x_{s})\,ds\big]\,dt\big|+
\sup_{\partial B_{1}}|u|.
$$
and to obtain \eqref{9.5.6} from \eqref{9.5.5}
it only remains to use the probabilistic version of the
Alexandrov estimate once more. The theorem is proved.

\begin{remark}
                                     \label{remark 4.29.3}
One of consequences of Theorem \ref{theorem 4.29.1}
is obtained when one takes $f$ to be the indicator function
of a Borel $G\subset B_{1}$.
Then estimate \eqref{4.29.5} says that
$|G|^{1/\gamma}$ is less than a constant $N$
times the average time spent by $x_{t}$ in $G$ before
exiting from $B_{1}$, where $|G|$ is the Lebesgue measure of $G$.

It turns out that even in such estimates
of the average time spent by $x_{t}$ in $G$ before
exiting from $B_{1}$
the constant $\gamma$ can be very small
when $\delta$ is small, so that there is no hope
to get \eqref{4.29.5} with large $\gamma$ for arbitrary $f$.

For instance, take $d=2$, $b=c=0$,
$$
\bar{a}^{ij}(x)=\delta^{ij}-\varepsilon\frac{x^{i}x^{j}}{|x|^{2}}
$$
for $x\ne0$ and $\bar{a}^{ij}(x)=\delta^{ij}$,
where $\varepsilon=1-\delta$. Then let $e_{1}$ be the first
basis vector and set
$$
a^{ij}(x)=\bar{a}^{ij}(x-e_{1}/2)
$$
Also let $G$ be the indicator of $B_{r}+e_{1}/2$, where $r\in(0,1/2)$.

Next solve the equation
\begin{equation}
                                               \label{4.29.9}
a^{ij}D_{ij}u(x)=-1
\end{equation}
in $B_{3/2}+e_{1}/2$ with zero boundary condition.
Then the value at zero of this solution will be
the average time spent by $x_{t}$ in $G$ before
exiting from $B_{3/2}+e_{1}/2$ and since the latter
contains $B_{1}\supset G$, $u(0)$ is greater than the average
time spent by $x_{t}$ in $G$ before leaving $B_{1}$.
It turns out that
\begin{equation}
                                               \label{4.29.3}
u(0)=\frac{1-\varepsilon}{\varepsilon(2-\varepsilon)}
2 ^{\varepsilon/(1-\varepsilon)}(
1-3 ^{-\varepsilon/(1-\varepsilon)})r^{(2-\varepsilon)/(1-\varepsilon)}
\end{equation}
which equals a constant times $|G|^{1/\gamma}$ with
$$
\gamma=\frac{2(1-\varepsilon)}{2-\varepsilon}.
$$
Thus, $\gamma$ can be made as small as we wish on the account
of taking $\delta$ small enough or $\varepsilon$ close to $1$.

One solves \eqref{4.29.9} in polar coordinates with pole
at $e_{1}/2$. Then, if $\rho$ is the polar radius, our
function $u(x)$ is written as $v(\rho)$ and $v$ satisfies
$$
(1-\varepsilon)v''+\frac{1}{\rho}v'=-I_{[0,r]}(\rho)
$$
with boundary conditions $v'(0)=0$ and $v(3/2)=0$.
The latter equation is easily solvable by using an
appropriate integrating factor, yields a function
$v$ with bounded second-order derivative, and after noting that 
$u(0)=v(1/2)$ leads to \eqref{4.29.3}.
\end{remark}

\mysection{Estimates in $L_{p}$  of  resolvent operators.
Parabolic case}

For a domain $Q\in\bR^{d+1}$ denote by $\partial'Q$
the parabolic boundary of $Q$, that is the set of
all points $(t_{0},x_{0})\in\partial Q$, for each of which
there exists a $\kappa>0$ and a continuous $\bR^{d}$-valued
function $x(t)$ defined on $[t_{0}-\kappa,t_{0}]$ such that
$x(t_{0})=x_{0}$ and $(t,x(t))\in Q$ for 
 $t\in[t_{0}-\kappa,t_{0})$. In case $Q=\bR^{d+1}$
we have $\partial Q=\partial'Q=\emptyset$.

Take $p\in[d+1,\infty)$  and  
introduce
$$
\hat W^{1,2}_{p}(Q)=\bigcap_{G\subset Q}W^{1,2}_{p}(G),
$$
where the intersection is taken over all bounded 
open subsets
$G$ of $Q$.

Set
$$
W^{1,2}_{p}= W^{1,2}_{p}(\bR^{d+1}),\quad
L_{p}=L_{p}(\bR^{d+1})
$$
and denote by $C(\bar{Q})$ the set of bounded continuous
functions on $\bar{Q}$.
Next, let
$$
L_{0}=\partial_{t}+a^{ij}(t,x)D_{ij},
$$
where $a(t,x)=(a^{ij}(t,x))$ is a $d\times d$ symmetric
matrix-valued function.
Let $\bR^{d}$-valued function $b(t,x)=(b^{1}(t,x),...,b^{d}(t,x))$ and
real-valued function $c(t,x)$ be defined
on $\bR^{d+1}$. Set
$$
L=L_{0} +b^{i}D_{i}-c,
$$
fix a $\delta>0$ and $K\in[0,\infty)$ and impose the following.
\begin{assumption}
                                               \label{assumption 8.29.1}
(i) For any $\xi\in\bR^{d}$ and all values of the arguments
$$
a^{ij}(t,x)\xi^{i}\xi^{j}\geq\delta|\xi|^{2},\quad
\tr a(t,x)+1\leq K,
$$

(ii) We have 
 $b=b_{1}+b_{2}$, where $b_{1}$ is bounded and $b_{2}\in L_{d+1}$.
 The function $c$ is nonnegative and bounded.

\end{assumption}

Our main goal in this section is to establish estimates
like
\begin{equation}
                                                       \label{4.24.1}
\lambda(\mu)\|u\|_{L_{p}(Q)}\leq N\|\mu u-Lu\|_{L_{p}(Q)},
\end{equation}
for $u\in W^{1,2}_{p}(Q)$ vanishing on $\partial'Q$,
where $N$ is a constant and the function $\lambda(\mu)>0$
for $\mu>0$ behaves like $\mu$ as $\mu\to\infty$.
The linear behavior of $\lambda(\mu)$ for large $\mu$  is, of course,
the best one could expect.

In a particular case of bounded $b\not\equiv0$ as we will see
from Corollary \ref{corollary 9.19.1} one can take
$\lambda(\mu)=\mu^{2}$ for $\mu$ close to $0$.

\begin{remark}
                                         \label{remark 9.12.3}
For $\theta\in(0,\infty)$ we introduce   $ 
\mu_{\theta}(\lambda) $ as a continuous 
nonnegative increasing function
of $\lambda>0$ such that
\begin{equation}
                                                  \label{4.27.1}
\|(|b| -\mu_{\theta}(\lambda))_{+}\|_{L_{d+1}}\leq \theta
\lambda^{-1/(2d+2)}.
\end{equation}

By our Assumption \ref{assumption 8.29.1} (ii), for any $\theta,
\lambda\in(0,\infty)$, there exists a
such a $\mu_{\theta}(\lambda)$ satisfying \eqref{4.27.1}.

On the other hand,
if, for some $\theta,\lambda\in(0,\infty)$, we can find an appropriate
constant $\mu$, then 
Assumption \ref{assumption 8.29.1} (ii) is satisfied
and one can find $\mu_{\theta}(\lambda)$ satisfying \eqref{4.27.1}
for all $\theta,\lambda\in(0,\infty)$.
Indeed,   then
$(|b|-\mu)_{+}\in L_{d+1}$ for some $\mu\in[0,\infty)$.
The latter means that $b=b_{1}+b_{2}$, where 
$$
b_{1}=bI_{|b|\leq\mu }+\frac{b}{|b|}\mu I_{|b|>\mu }
$$ 
is bounded
and 
$$
b_{2}=bI_{|b|>\mu }-\frac{b}{|b|}\mu I_{|b|>\mu }  \in L_{d+1}.
$$

To satisfy our requirement for $\mu_{\theta}(\lambda)$
to be increasing and continuous, as is easy to see,
 one can just take
$$
\mu_{\theta}(\lambda)
=\inf\{\mu\geq0:\|(|b|   
-\mu)_{+}\|_{L_{p}}\leq \theta\lambda^{-1/(2d+2)}\} .
$$

\end{remark}

 \begin{remark}
                                         \label{remark 4.28.1}
If $b_{2}\in L_{d+2}$, then for any $\nu>0$
$$
\int_{\bR^{d+1}}\nu(|b_{2}|-\nu)_{+}^{d+1}\,dxdt
\leq
\int_{\bR^{d+1}}  |b_{2}| ^{d+2}\,dxdt.
$$
In particular, for any $\lambda,\nu\in(0,\infty)$ and $M:=\sup|b_{1}|$
$$
\int_{\bR^{d+1}} (|b|-\nu \lambda^{1/2}-M)_{+}^{d+1}\,dxdt
\leq (\|b_{2}\|_{L_{d+2}}^{d+2}/\nu )\lambda^{-1/2}
$$
and one can take $\nu_{\theta}\lambda^{1/2}+M$ as $\mu_{\theta}
(\lambda)$ in Remark \ref{remark 9.12.3} if one chooses
$$
\nu_{\theta}\geq\|b_{2}\|_{L_{d+2}}^{d+2}\theta^{-d-1}.
$$

\end{remark}

In the following 
main result of the section the case $Q=\bR^{d+1}$ is allowed.
Of course, in that case no conditions on the values of $u$
on $\partial' Q$ are necessary.
 
\begin{theorem}
                                               \label{theorem 8.25.1}
There is a constant $\theta=\theta(\delta,d)>0$ such that,
if $\lambda>0$, $u\in \hat W^{1,2}_{d+1 }(Q)\cap C(\bar{Q})$,  
$u\leq 0$ on $\partial'Q$,
and in case $Q$ is unbounded
\begin{equation}
                                                      \label{5.1.1}
\nlimsup_{\substack{t+|x|\to\infty,\\(t,x)\in Q}}u (t,x)\leq 0,
\end{equation}
then we have:

(i) For $p=d+1$ and $\mu\geq
K\lambda+\mu_{\theta}(\lambda)\lambda^{1/2}$, it holds that
\begin{equation}
                                                      \label{8.29.5}
\lambda\|u_{+}\|_{L_{p}(Q)}\leq N\|(\mu u-Lu)_{+}\|_{L_{p}(Q)};
\end{equation}

(ii) For $p\geq d+1$ and $\mu\geq
K\lambda+\mu_{\theta}(\lambda)\lambda^{1/2}$, it holds that
\begin{equation}
                                                      \label{4.28.2}
\lambda^{(d+1)/p}
\mu^{1-(d+1)/p}\|u_{+}\|_{L_{p}(Q)}
\leq N\|(\mu u-Lu)_{+}\|_{L_{p}(Q)}.
\end{equation}

  In all cases $N=N(d,p,\delta )$.
\end{theorem} 

Before we prove this theorem we extract a few corollaries.

\begin{corollary}
                                             \label{corollary 5.1.3}
If $b_{2}\in L_{d+2}$, then for any 
$p\geq d+1$, $\mu>0$, and $u\in\hat  W^{1,2}_{d+1} (Q)
\cap C(\bar{Q})$,
such that  
$u\leq0$ on $\partial'Q$ and condition  
 \eqref{5.1.1} is satisfied,   
 we have
$$
\mu\|u_{+}\|_{L_{p}(Q)}\leq N(K+\nu_{\theta}
+1)^{(d+1)/p}\|(\mu u-Lu)_{+}\|_{L_{p}(Q)} 
$$
for  $\mu\geq(K+\nu_{\theta}+1)M^{2} $,
where $\nu_{\theta}$ and $M$ are 
taken from Remark \ref{remark 4.28.1}, and
$$
\mu^{1+(d+1)/p}\|u_{+}\|_{L_{p}(Q)}\leq N
[(K+\nu_{\theta}+1)M]^{2(d+1)/p}\|(\mu u-Lu)_{+}\|_{L_{p}(Q)}
$$
for  $\mu\leq(K+\nu_{\theta}+1)M^{2}$.
In all cases $N=N(d,p,\delta )$.
\end{corollary}

Indeed, take $\mu_{\theta}(\lambda)$ from Remark \ref{remark 4.28.1}.
Then for any $\mu>0$ one can find $\lambda(\mu)>0$ such that
$$
\mu= K\lambda(\mu)+\mu_{\theta}(\lambda)\sqrt{\lambda(\mu)}
=(K+\nu_{\theta})\lambda(\mu)+M\sqrt{\lambda(\mu)},
$$
which implies that \eqref{4.28.2} holds with $\lambda=\lambda(\mu)$.
After that it only remains to prove that 
$$
\lambda(\mu)\geq \mu/(K+\nu_{\theta}+1)\quad\text{if} \quad
\mu\geq(K+\nu_{\theta}+1) M^{2},
$$
$$\lambda(\mu)
\geq\mu^{2}(K+\nu_{\theta}+1)^{-2}M^{-2}\quad\text{if} \quad
0<\mu\leq (K+\nu_{\theta}+1)M^{2}.
$$
 This is 
easily done after observing that
$x:=\mu/M^{2}$ and $y:=\sqrt{\lambda(\mu)}/M$ satisfy
$x=(K+\nu_{\theta})y^{2}+y$
  and the rest is left to the reader.

Here is a particular case of Corollary \ref{corollary 5.1.3}
when $p=d+1$.

\begin{corollary}
                                           \label{corollary 9.17.03}
If $b_{2}\in L_{d+2}$, then for any
$\mu>0$ and $u\in\hat  W^{1,2}_{d+1} (Q)
\cap C(\bar{Q})$,
such that  
$u\leq0$ on $\partial'Q$ and condition  
 \eqref{5.1.1} is satisfied,  
 we have
$$
\mu\|u_{+}\|_{L_{d+1}(Q)}\leq N(K+\nu_{\theta}+1)\|(\mu
u-Lu)_{+}\|_{L_{d+1}(Q)}  
$$
if $\mu\geq (K+\nu_{\theta}+1)M^{2}$, and
$$
\mu^{2}\|u_{+}\|_{L_{d+1}(Q)}\leq N(K+\nu_{\theta}+1)^{2}
M^{2}\|(\mu
u-Lu)_{+}\|_{L_{d+1}(Q)}   
$$
if $ \mu\leq (K+\nu_{\theta}+1)M^{2}$,
where $N=N(d,\delta )$.
\end{corollary}

In the case of bounded $b$ we have a version
of Theorem \ref{theorem 8.25.1}, which is easier to memorize.
The first part of
the following result was used, for instance, in \cite{DKL}.

\begin{corollary}
                                      \label{corollary 9.19.1}
  Assume that $b$ is bounded and set $M=\sup|b|$.
Then for any $\mu>0$ and 
$u\in \hat W^{1,2}_{d+1}(Q)\cap C(\bar{Q})$,
such that  
$u\leq 0$ 
on $\partial'Q$ and condition  \eqref{5.1.1} 
 is  satisfied,   we have
$$
\mu\|u_{+}\|_{L_{p}(Q)}\leq N(K+1)\|(\mu
u-Lu)_{+}\|_{L_{p}(Q)} \quad\text{if}
\quad \mu\geq (K+1)M^{2},
$$
$$
\mu^{2}\|u_{+}\|_{L_{p}(Q)}\leq N(K+1)^{2}M^{2}\|(\mu
u-Lu)_{+}\|_{L_{p}(Q)} \quad\text{if}
\quad \mu\leq (K+1)M^{2},
$$
where $N=N(\delta,d,K)$.
\end{corollary}

 Indeed, we have $\nu_{\theta}=0$  and
$(K+1)^{(d+1)/p}\leq K+1$ whereas
$$
\mu^{1+(d+1)/p}\geq \mu^{2}[(K+1)M]^{-2+2(d+1)/p}
$$
if  $\mu\leq(K +1)M^{2}$.
 
\begin{remark}
The fact that in Corollary
\ref{corollary 9.19.1}  we have the factors $\mu$ and $\mu^{2}$
in different ranges of $\mu$ may look suspicious.
In  Example \ref{example 9.20.1} we give an argument
partially explaining this effect.
\end{remark}

In the general case, if $b_{2}\in L_{d+1}$ one can still
obtain estimates like in Corollaries \ref{corollary 9.17.03}
 and \ref{corollary 9.19.1}
for $\mu$ small.

\begin{corollary}
                                           \label{corollary 9.17.2}
Take a $\lambda'>0$ and set
$\mu':=K\lambda'+\mu_{\theta}(\lambda')
\sqrt{\lambda'}$. Then for any  $0<\mu\leq\mu'$ and
  $u\in \hat W^{1,2}_{d+1}(Q)\cap C(Q)$,
such that $u\leq 0$ 
on $\partial'Q$ and condition \eqref{5.1.1} 
 is  satisfied,  
we have
\begin{equation}
                                               \label{9.17.3}
\mu^{1+(d+1)/p} \|u_{+}\|_{L_{p}(Q)}\leq N
(\mu')^{ 2(d+1)/p}(\lambda')^{-(d+1)/p}
\|(\mu u-Lu)_{+}\|_{L_{p}(Q)}, 
\end{equation}
where $N$ is the constant in \eqref{4.28.2}.
 \end{corollary}

Indeed, define $\lambda>0$ from 
\begin{equation}
                                                       \label{4.29.1}
\mu= K\lambda+\mu_{\theta}(\lambda)
\sqrt{\lambda}.
\end{equation}
 This is possible since $\mu_{\theta}(\lambda)$
is an increasing and continuous function of $\lambda$.
Then \eqref{4.28.2} holds. Since $\mu\leq\mu'$ we have that
$\lambda\leq\lambda'$ and
$$
\mu= K\lambda+\mu_{\theta}(\lambda)\sqrt{\lambda}\leq
\sqrt{\lambda}(K\sqrt{\lambda'}+\mu_{\theta}(\lambda'))
=\sqrt{\lambda}\mu'(\lambda')^{-1/2},
$$
$$
\lambda\geq \mu^{2}(\mu')^{-2}\lambda'.
$$
Finally,
$$
\lambda^{(d+1)/p}
\mu^{1-(d+1)/p}\geq \mu^{1+(d+1)/p}
(\mu')^{-2(d+1)/p}(\lambda')^{(d+1)/p}.
$$

\begin{remark}
                                                \label{remark 4.29.1}

Unfortunately, in general, there is no control on how
fast $\mu_{\theta}(\lambda)$ may grow to infinity
as $\lambda\to\infty$. Accordingly, we do not
know how slow the solution of
\eqref{4.29.1} may go to infinity as $\mu\to\infty$,
so that we were able to prove the natural rate $\mu^{-1}$ of decay
of the resolvent operator $R_{\mu}$ of $L$ only if $b_{2}\in L_{d+2}$.
Actually, the author 
conjectures that for some $b$, if $p\in[p+1,p+2)$,
the $L_{p}\to L_{p}$-norm of $R_{\mu}$ may have as slow
power decay
   as we wish as $\mu\to\infty$. Still 
from our results we have that
$\lambda\to\infty$ as $\mu\to\infty$ so that the norm
of $R_{\mu}$ as an operator in $L_{p}$ does go to zero
as $\mu\to\infty$.
We will see later that  
in the elliptic case we will not have this issue.

\end{remark}

{\bf Proof of Theorem \ref{theorem 8.25.1}}. 
Take a $\varepsilon>0$ and define 
$$
Q_{\varepsilon}
=\{(t,x)\in Q:u(t,x)>\varepsilon\}.
$$
Obviously $Q_{\varepsilon}$ is a bounded domain,
 $u-\varepsilon=0$ on $\partial'Q$,
and $u-\varepsilon\in W^{1,2}_{d+1}(Q_{\varepsilon})$.
If the assertions of the theorem are true when $Q$ is bounded,
$u\in  W^{1,2}_{d+1}(Q)\cap C(\bar{Q})$, and $u\leq0$
on $\partial'Q$, then, applying them to
$Q_{\varepsilon}$ and $u-\varepsilon$ and passing to the limit
as $\varepsilon\downarrow0$ 
on the basis of the monotone convergence theorem,
we obtain the assertions in full generality.
Therefore, we may assume that $Q$ is bounded,
$u\in  W^{1,2}_{d+1}(Q)\cap C(\bar{Q})$, and $u\leq0$
on $\partial'Q$.

Next let
 $b_{n}=bI_{|b|\leq n}$ and    observe that the original
$\mu_{\theta}(\lambda)$ is satisfying \eqref{4.27.1} with $b_{n}$ 
in place of
$b$.
Hence, if the theorem is true under one more
 additional assumption that
$b$ is bounded, then under the conditions in (i)
\begin{equation}
                                                      \label{9.13.2}
\lambda\|u_{+}\|_{L_{d+1}(Q)}\leq N(d, \delta )
\|(\mu u-L^{n}u)_{+}\|_{L_{d+1}(Q)},
\end{equation}
where $L^{n}=L+(b_{n}^{i}-b^{i})D_{i}$.
Observe that if $(b^{i}D_{i}u)_{-}\not\in L_{d+1}(Q)$,
then the right-hand side of \eqref{8.29.5} is infinite
and we have nothing to prove. In case $(b^{i}D_{i}u)_{-} \in L_{d+1}(Q)$,
  we can pass to the limit 
in \eqref{9.13.2} and obtain \eqref{8.29.5}.
Similar situation occurs in the case of   assertion (ii)
of the theorem.
Therefore, in the rest of the proof of the theorem
  we may assume
that $b$  is bounded.

Using approximations we convince
ourselves that we may also assume that   $a^{ij},c\in
C^{\infty}_{b}$. In that case
introduce $I$ as the set of $\mu>0$ for each of which
the operator
$\mu-L$ as an operator from $W^{1,2}_{d+1}$
to $L_{d+1}$ is onto, invertible, and the inverse 
$R_{\mu}:=(\mu-L)^{-1}$ is bounded as an operator from
$L_{d+1}$ onto $W^{1,2}_{d+1}$. Obviously $I$ is an open subset 
of $(0,\infty)$.
It is well known (see, for instance, \cite{Kr08})
that  all large $\mu$ are in $I$. Therefore,
it makes sense to introduce $\mu'$ as the smallest
number such that $(\mu',\infty)\subset I$.

 Also notice that,
if $u\in W^{1,2}_{d+1}(Q)\cap C(\bar{Q})$,
$u\leq 0$ on $\partial'Q$, 
and $\mu>\mu'$, 
then by the maximum principle
(see, for instance, Theorem 3.4.2 of \cite{Kr87}) in $Q$ we have
$$
u\leq R_{\mu}[(\mu u-Lu)_{+}I_{Q}],
\quad u_{+}\leq R_{\mu}[(\mu u-Lu)_{+}I_{Q}].
$$
It follows that for any $p$
$$
\|u_{+}\|_{L_{p}(Q)}\leq\|R_{\mu}[(\mu u-Lu)_{+}I_{Q}]\|_{L_{p}}
$$
and reduces the proof of the theorem
to proving that $\mu'< K\lambda+\mu_{\theta}(\lambda)
\lambda^{1/2}$ for any $\lambda>0$ and that
\begin{equation}
                                         \label{5.1.4}
\lambda^{(d+1)/p}
\mu^{1-(d+1)/p}\|R_{\mu}f\|_{L_{p}}\leq N\|f\|_{L_{p}}
\end{equation}
as long as $f\geq0$ and $\mu\geq
K\lambda+\mu_{\theta}(\lambda)\lambda^{1/2}$.

First we deal with $p=d+1$ and then we use the
 Marcinkiewicz interpolation theorem.
For $\mu>\mu'$
denote by $N_{\mu}$ the norm of  $R_{\mu}$
as an operator acting from $L_{d+1}$ {\em into\/} $L_{d+1}$,
that is the least constant $N$ such that
\begin{equation}
                                                      \label{9.5.8}
\|R_{\mu}g\|_{L_{d+1}}\leq N\|g\|_{L_{d+1}}
\end{equation}
for all $g\in L_{d+1}$. 

Our first goal is to show that
\begin{equation}
                                                      \label{9.10.5}
N_{\mu}\leq N\lambda^{-1}
\end{equation}
as long as $\mu>\mu'$ and $\mu 
\geq
K\lambda+\mu_{\theta}(\lambda)\lambda^{1/2}$, where $N=N(d,\delta )$,
  $\lambda>0$ is fixed, and $\theta(\delta,d)>0$ is to be chosen
appropriately.

By the maximum principle
$|R_{\mu}g|\leq R_{\mu}|g|$, so that $N_{\mu}$ 
is also the least constant $N$ for which
\eqref{9.5.8} holds for all {\em
nonnegative\/}  $g\in L_{d+1}$.

Observe that owing to \cite{Kr74} 
 for any   nonnegative $f\in C^{\infty}_{0}$ 
there exists a nonnegative function $\psi_{\lambda}$, which is
$\lambda$-convex in $x$, decreasing in $t$ and such that
for any  $\varepsilon>0$ (see equation (29) in \cite{Kr74})
\begin{equation}
                                                      \label{8.26.1}
L_{0}\psi^{\varepsilon}_{\lambda}-\lambda(\tr a+1) 
\psi^{\varepsilon}_{\lambda}\leq -f^{\varepsilon},
\end{equation}
where the notation $v^{\varepsilon}$ stands for a standard
mollification of $v$ with kernel of support diameter $\varepsilon$.
 Furthermore,
\begin{equation}
                                                      \label{8.29.1}
\sup \psi_{\lambda}\leq N\lambda^{-d/(2d+2)}\|f\|_{L_{d+1}},
\end{equation}
\begin{equation}
                                                      \label{8.29.2}
\|\psi_{\lambda}\|_{L_{d+1}}\leq N\lambda^{-1}\|f\|_{L_{d+1}},
\end{equation}
where the first estimate is a combination of
estimates (12) and (13)  of \cite{Kr74} and the second one
is obtained before Theorem 3  of \cite{Kr74}. The above cited results
of \cite{Kr74} are obtained there by using the  theory of controlled
diffusion processes. The more PDE oriented reader may like to use
Theorem 3.2.8 of \cite{Kr87}. We emphasize that the above constants
$N$ depend only on $d$ and $\delta$.

We also know that $|D\psi_{\lambda}^{\varepsilon}|\leq
\psi_{\lambda}^{\varepsilon}\sqrt{\lambda}$.
Therefore, \eqref{8.26.1}  implies that
\begin{equation}
                                                      \label{8.26.2}
L\psi^{\varepsilon}_{\lambda}-\mu
\psi^{\varepsilon}_{\lambda}\leq -f^{\varepsilon}
+(|b|\sqrt{\lambda}+\lambda(\tr a+1) -\mu)\psi^{\varepsilon}_{\lambda}.
\end{equation}

By the maximum principle
$$
\psi_{\lambda}\geq R_{\mu}f^{\varepsilon}-
R_{\mu}(|b|\sqrt{\lambda}+\lambda(\tr a+1)
-\mu)\psi^{\varepsilon}_{\lambda}
$$
and \eqref{8.29.1} and \eqref{8.29.2} yield
$$
 \|R_{\mu}f^{\varepsilon}\|_{L_{d+1}}
\leq N\lambda^{-1}\|f\|_{L_{d+1}}
$$
$$
+
N\lambda^{-d/(2d+2)}\|f\|_{L_{d+1}}
\|R_{\mu}(|b|\sqrt{\lambda}+\lambda(\tr a+1) -\mu)_{+}\|_{L_{d+1}}
$$
\begin{equation}
                                                      \label{8.29.3}
\leq N'\|f\|_{L_{d+1}}\big(\lambda^{-1}
+\lambda^{-d/(2d+2)}N_{\mu}\|(|b|\sqrt{\lambda}+\lambda(\tr a+1) 
-\mu)_{+}\|_{L_{d+1}}\big),
\end{equation}
where $N'=N'(d,\delta  )$.

By using the Alexandrov estimate and Fatou's lemma we can pass to the
limit in \eqref{8.29.3} as $\varepsilon\downarrow0$ and then
we obtain    
$$
 \|R_{\mu}f\|_{L_{d+1}}
\leq   N'\|f\|_{L_{d+1}}\big(\lambda^{-1}
$$
\begin{equation}
                                                      \label{9.5.9}
+\lambda^{-d/(2d+2)}N_{\mu}\|(|b|\sqrt{\lambda}+\lambda(\tr a+1) 
-\mu)_{+}\|_{L_{d+1}}\big).
\end{equation}
We have proved \eqref{9.5.9} for nonnegative $f\in C^{\infty}_{0}$.
The Alexandrov estimate and Fatou's lemma allow us to carry over
this estimate to arbitrary nonnegative $f\in L_{d+1}$.
After that we recall what was said about $N_{\mu}$ in connection
with \eqref{9.5.8} and by the definition of $N_{\mu}$  we conclude that
\begin{equation}
                                                      \label{8.29.4}
N_{\mu}\leq 
N'\lambda^{-1}
+N'\lambda^{-d/(2d+2)}N_{\mu}\|(|b|\sqrt{\lambda}+K\lambda  
-\mu)_{+}\|_{L_{d+1}}.
\end{equation}
For $\mu\geq
K\lambda+\mu_{\theta}(\lambda)\lambda^{1/2}$
the factor of $N_{\mu}$ in \eqref{8.29.4}  is dominated by
$$
N'\lambda^{-d/(2d+2)}\|(|b|\sqrt{\lambda}  
-\mu_{\theta}(\lambda)\sqrt{\lambda})_{+}\|_{L_{d+1}}
$$
$$
=N'\lambda^{1/(2d+2)}\|(|b| 
-\mu_{\theta}(\lambda))_{+}\|_{L_{d+1}},
$$
which by the definition of $\mu_{\theta}(\lambda)$
 is less than or equal to
$N'\theta$. Thus, if $\theta=\theta(d,\delta )>0$ is chosen in such a way that
$N'\theta\leq1/2$, then $N_{\mu}\leq 2N'\lambda^{-1}$,
which is \eqref{9.10.5}. Thus, \eqref{9.10.5} holds
if $\mu>\mu'$ and $\mu\geq
K\lambda+\mu_{\theta}(\lambda)\lambda^{1/2}$.

To finish considering the case that
$p=d+1$ it suffices to show that 
$$
\mu'
< K\lambda+\mu_{\theta}(\lambda)\lambda^{1/2}.
$$ 
To this end suppose that $\mu'\geq 
K\lambda+\mu_{\theta}(\lambda)\lambda^{1/2}$.
Then by the above for any $u\in W^{1,2}_{d+1}$ the inequality
\begin{equation}
                                                      \label{9.13.3}
\lambda\|u \|_{L_{d+1}}\leq N\| \mu u-Lu \|_{L_{d+1}}
\end{equation}
holds
if $\mu>\mu'$ and by continuity if $\mu=\mu'$
as well. Furthermore, as is well known, under our additional assumptions
on the coefficients of $L$, there are constants $M_{i}<\infty$ such that
for any $u\in W^{1,2}_{d+1}$
$$
\|u\|_{W^{1,2}_{d+1}}\leq M_{1}(\|Lu\|_{L_{d+1}}+\|u\|_{L_{d+1}})
 \leq M_{2}(\|\mu u-Lu\|_{L_{d+1}}+(1+\mu)\|u\|_{L_{d+1}}).
$$
Owing to \eqref{9.13.3} (recall that $\lambda>0$ is fixed)
$$
\|u\|_{W^{1,2}_{d+1}}\leq M_{1}(\|Lu\|_{L_{d+1}}+\|u\|_{L_{d+1}})
 \leq M_{3}(1+\mu) \|\mu u-Lu\|_{L_{d+1}},
$$
where $M_{3}$ is independent of $u$ and $\mu\geq\mu'$. By the method of
continuity applied with respect to $\mu$ this estimate
implies that $\mu'\in I$, which yields the desired
contradiction with the definition
of $\mu'$. This proves that \eqref{5.1.4} holds for $p=d+1$
and any $f\in L_{p}$ as long as $\mu\geq
K\lambda+\mu_{\theta}(\lambda)\lambda^{1/2}$.

Next observe that the maximum principle
implies that   $R_{\mu}$ is well defined as an operator
in $L_{\infty}$ and its norm is less than or equal to $\mu^{-1}$
for any $\mu>0$. Then 
we obtain that \eqref{5.1.4} holds for $p\geq d+1$
and all $f\in L_{p}$ by the 
 Marcinkiewicz interpolation theorem.
The positivity of the operator $R_{\mu}$ and the monotone
convergence theorem allows us to conclude that
\eqref{5.1.4} holds for $p\geq d+1$
and all $f\geq0$ as long as $\mu\geq
K\lambda+\mu_{\theta}(\lambda)\lambda^{1/2}$. 
The theorem is proved.

\mysection{Estimates in $L_{p}$  of the resolvent operators.
Elliptic case}

Take $p\in[d ,\infty)$, a domain $Q\subset\bR^{d}$,
introduce $W^{2}_{p}(Q)$ as the usual Sobolev space
and $\hat W^{2}_{p}(Q)$ as the collection of all $u$
which belong to $W^{2}_{p}(G)$ for any bounded subdomain of $G$.
 Also denote
$$
W^{ 2}_{p}= W^{ 2}_{p}(\bR^{d }),\quad
L_{p}=L_{p}(\bR^{d }).
$$
Introduce
$$
L_{0}= a^{ij}( x)D_{ij},
$$
where $a( x)=(a^{ij}( x))$ is a $d\times d$ symmetric
matrix-valued function.
Let $\bR^{d}$-valued function $b( x)=(b^{1}( x),...,b^{d}( x))$ and
real-valued function $c( x)$ be defined
on $\bR^{d }$. Set
$$
L=L_{0} +b^{i}D_{i}-c,
$$
fix a $\delta>0$ and $K\in[0,\infty)$ and impose the following.
\begin{assumption}
                                               \label{assumption 9.10.1}
(i) For any $x,\xi\in\bR^{d}$ 
$$
a^{ij}( x)\xi^{i}\xi^{j}\geq\delta|\xi|^{2},\quad
\tr a( x) \leq K,
$$

(ii) We have $b=b_{1}+b_{2}$, where $b_{1}$
is bounded and $b_{2}\in L_{d}$.

 The function $c$ is nonnegative and bounded.

\end{assumption}

\begin{remark}
                                            \label{remark 9.15.1}

For $\theta\in(0,\infty)$ introduce $\mu_{\theta}\in(0,\infty)$
in such a way that
$$  
 \|(|b|   
-\mu_{\theta})_{+}\|_{L_{d} }\leq \theta .
$$
 
As in Remark \ref{remark 9.12.3} Assumption \ref{assumption 9.10.1}(ii)
is satisfied if and only if 
there exists a $\theta\in(0,\infty)$ such that
one can find an appropriate $\mu_{\theta}$ and in this case
one can find an appropriate $\mu_{\theta}$
 for any $\theta\in(0,\infty)$. 
\end{remark}

In the following theorem the case $Q=\bR^{d}$ is allowed.
Of course, in that case no conditions on the values of $u$
on $\partial Q$ are necessary.
\begin{theorem}
                                               \label{theorem 9.10.1}
There exists a constant $\theta=\theta(\delta,d )>0$
such that, if $u\in \hat W^{ 2}_{d}(Q)\cap C(\bar{Q})$,
$u\leq 0$ on $\partial Q$,
and in case $Q$ is unbounded 
\begin{equation}
                                                      \label{5.1.2}
\nlimsup_{\substack{ |x|\to\infty,\\x\in Q}}u (t,x)\leq 0,
\end{equation}
then for any
 $\lambda>0$ and 
  $\mu\geq
K\lambda+\mu_{\theta}\lambda^{1/2}$  we have  
\begin{equation}
                                                      \label{9.10.6}
\lambda^{d/p}
\mu^{1-d/p}\|u_{+}\|_{L_{p}(Q)}
\leq N\|(\mu u-Lu)_{+}\|_{L_{p}(Q)},  
\end{equation}
where $N=N(d,\delta,p )$.
\end{theorem}

From the proof of Corollary
\ref{corollary 5.1.3} we know that
if $\mu=
K\lambda+\mu_{\theta}\lambda^{1/2}$, then
$$
\lambda\geq\mu/(K+1)\quad\text{for}\quad \mu\geq(K+1)\mu^{2}_{\theta},
$$
$$
\lambda\geq\mu^{2} (K+1)^{-2}\mu_{\theta}^{-2}
\quad\text{for}\quad \mu\leq(K+1)\mu^{2}_{\theta}.
$$
This leads to the following.  

\begin{corollary}    
                                               \label{corollary 9.15.1}
For any $\mu>0$ and 
 $u\in \hat W^{ 2}_{d}(Q)\cap C(\bar{Q})$, such that
$u\leq 0$ on $\partial Q$ and \eqref{5.1.2} holds, we have
$$
\mu\|u_{+}\|_{L_{p}(Q)}
\leq N (K+1)^{d/p}\|(\mu u-Lu)_{+}\|_{L_{p}(Q)}
$$
if $\mu\geq(K+1)\mu_{\theta}^{2}$ and
$$
\mu^{1+d/p}\|u_{+}\|_{L_{p}(Q)}\leq 
N [(K+1)\mu_{\theta}]^{2d/p}\|(\mu u-Lu)_{+}\|_{L_{p}(Q)}
$$
if $\mu\leq(K+1)\mu_{\theta}^{2}$,
where $N$ is the constant in \eqref{9.10.6}.

\end{corollary}

In the case of bounded $b$ our results lead to 
 a simpler statement.

\begin{corollary}
                                                  \label{corollary 9.20.3}
  Assume that $b$ is bounded and set $M=\sup|b|$.
Then for any 
$\mu>0$ and 
$u\in\hat W^{ 2}_{d}(Q)\cap C(\bar{Q})$, such that
$u\leq 0$ on $\partial Q$ and \eqref{5.1.2} holds,  we have
$$
\mu\|u_{+}\|_{L_{p}(Q)}\leq N(K+1)\|(\mu
u-Lu)_{+}\|_{L_{p}(Q)} \quad\text{if}
\quad \mu\geq (K+1)M^{2},
$$
$$
\mu^{2}\|u_{+}\|_{L_{p}(Q)}\leq N(K+1)^{2}M^{2}\|(\mu
u-Lu)_{+}\|_{L_{p}(Q)} \quad\text{if}
\quad \mu\leq (K+1)M^{2},  
$$
where $N$ is the constant in \eqref{9.10.6}.
\end{corollary}

The proof is
almost identical to the proof of Corollary \ref{corollary 9.19.1}.
It turns out that for $p=d=1$ the result of Corollary
\ref{corollary 9.20.3} (especially concerning the case
of small $\mu$) is rather sharp.  

\begin{example}
                                               \label{example 9.20.1}
Take $d=p=1$
a constant $M>0$ and for $\mu>0$ consider the equation
$u''+bu'-\mu u=-f$, where $b(x)=-M\text{\rm sign\,} x$.
If $f$ approaches the $\delta$-function concentrated at the origin,
its $L_{1}$-norm tends to one. The limit of $L_{1}$-norms of the 
corresponding solution will be the $L_{1}$-norm of the fundamental
solution with ``pole'' at the origin. This solution
is $e^{-\nu |x|}/(2\nu )$, where
$$
\nu =\frac{\sqrt{M^{2}+4\mu}-M}{2}=\frac{2\mu}{\sqrt{M^{2}+4\mu}+M}.
$$
Hence the limit $L_{1}$-norm of the solutions is
$$
\frac{1}{\nu ^{2}}=\frac{[\sqrt{M^{2}+4\mu}+M]^{2}}{4\mu^{2}}.
$$
The last expression should be multiplied by $\mu$ in order
to become bounded for large $\mu$ and by $\mu^{2}$ 
in order
to become bounded for small $\mu$. This shows that
the dichotomy is unavoidable and has exactly the form as in
Corollary \ref{corollary 9.20.3}.
\end{example}

{\bf Proof of Theorem \ref{theorem 9.10.1}}. 
We first concentrate on the case that
$p=d$.
As in the proof of Theorem \ref{theorem 8.25.1}
we may assume that $Q$ is bounded $u\in W^{2}_{p}(Q)\cap C(\bar{Q})$,
  $b$ is bounded,  and $a^{ij}$ and $c$ are infinitely differentiable
with bounded derivatives.
 In that case
it is well known (see, for instance, Theorem 11.6.2 in \cite{Kr08})
that for any $\mu>0$ the operator
$\mu-L$ as an operator from $W^{ 2}_{d}$
to $L_{d}$ is onto, invertible, and the inverse 
$R_{\mu}:=(\mu-L)^{-1}$ is bounded as an operator from
$L_{d}$ onto $W^{ 2}_{d}$.
Denote by $N_{\mu}$ the norm of  $R_{\mu}$
as an operator acting from $L_{d}$ {\em into\/} $L_{d}$,
that is the least constant $N$ such that
\begin{equation}
                                                      \label{9.10.3}
\|R_{\mu}g\|_{L_{d}}\leq N\|g\|_{L_{d}}
\end{equation}
for all $g\in L_{d}$. By the 
same reasons as in the proof of Theorem \ref{theorem 8.25.1} 
 to prove the theorem for $p=d$ we need only show that  
$N_{\mu}\leq N\lambda^{-1}$.

Observe that, if  $\lambda=1$, then, owing to \cite{Kr73},  
 for any    nonnegative $f\in C^{\infty}_{0}$ 
there exists a nonnegative function $\psi_{\lambda}$, which is
$\lambda$-convex in $x$ and 
\begin{equation}
                                                      \label{8.26.01}
L_{0}\psi^{\varepsilon}_{\lambda} - \lambda \tr a   
\psi^{\varepsilon}_{\lambda} \leq -f^{\varepsilon},
\end{equation}
where the notation $v^{\varepsilon}$ stands for a standard
mollification of $v$ with kernel of support diameter $\varepsilon$
(see the proof of Lemma 1 of \cite{Kr73}).
 Furthermore (see equation (22) in \cite{Kr73} and the end of the proof 
of Theorem 2 of \cite{Kr73}),
\begin{equation}
                                                      \label{8.29.01}
\sup \psi_{\lambda} \leq N\lambda^{-1/2} \|f\|_{L_{d}},
\end{equation}
\begin{equation}
                                                      \label{8.29.02}
\|\psi_{\lambda} \|_{L_{d}}\leq N\lambda^{-1} \|f\|_{L_{d}},
\end{equation}
where $N=N(d,\delta )$. Actually, dilations show that
one can take any $\lambda>0$. These results are obtained in \cite{Kr73}
by probabilistic methods. In terms of PDEs  
the existence of $\psi_{\lambda}$ with the properties described above
can be found in Theorem 3.2.3 of \cite{Kr87}.

We also know that $|D\psi_{\lambda}^{\varepsilon}|\leq
\psi_{\lambda}^{\varepsilon}\sqrt{\lambda}$.
Therefore \eqref{8.26.01}  implies that
\begin{equation}
                                                      \label{8.26.02}
L\psi^{\varepsilon}_{\lambda}-\mu
\psi^{\varepsilon}_{\lambda}\leq -f^{\varepsilon}
+(|b|\sqrt{\lambda}+\lambda \tr a  -\mu)\psi^{\varepsilon}_{\lambda}.
\end{equation}

By the maximum principle
$$
\psi_{\lambda}\geq R_{\mu}f^{\varepsilon}-
R_{\mu}(|b|\sqrt{\lambda}+\lambda \tr a 
-\mu)\psi^{\varepsilon}_{\lambda}
$$
and \eqref{8.29.01} and \eqref{8.29.02} yield
$$
 \|R_{\mu}f^{\varepsilon}\|_{L_{d}}
\leq N\lambda^{-1}\|f\|_{L_{d}}
$$
$$
+
N\lambda^{-1/2}\|f\|_{L_{d}}
\|R_{\mu}(|b|\sqrt{\lambda}+\lambda \tr a  -\mu)_{+}\|_{L_{d}}
$$
\begin{equation}
                                                      \label{8.29.03}
\leq N'\|f\|_{L_{d}}\big(\lambda^{-1}
+\lambda^{-1/2}N_{\mu}\|(|b|\sqrt{\lambda}+\lambda \tr a  
-\mu)_{+}\|_{L_{d}}\big),
\end{equation}
where $N'=N'(d,\delta )$.

By using the Alexandrov estimate and Fatou's lemma we can pass to the
limit in \eqref{8.29.03} as $\varepsilon\downarrow0$ and then
we obtain
$$
 \|R_{\mu}f\|_{L_{d}}
\leq N\lambda^{-1}\|f\|_{L_{d}}
$$
\begin{equation}
                                                      \label{9.5.09}
\leq N'\|f\|_{L_{d}}\big(\lambda^{-1}
+\lambda^{-1/2}N_{\mu}\|(|b|\sqrt{\lambda}+\lambda \tr a  
-\mu)_{+}\|_{L_{d}}\big).
\end{equation}

We have proved \eqref{9.5.09} for nonnegative $f\in C^{\infty}_{0}$.
The Alexandrov estimate and Fatou's lemma allow us to carry over
this estimate to arbitrary nonnegative $f\in L_{d}$.
As in the proof of Theorem \ref{theorem 8.25.1} constant $N_{\mu}$
is also the smallest constant for which \eqref{9.10.3} holds for all 
nonnegative $g$. Therefore, now
  \eqref{9.5.09} implies that
\begin{equation}  
                                                      \label{8.29.04}
N_{\mu}\leq 
N'\lambda^{-1}
+N'\lambda^{-1/2}N_{\mu}\|(|b|\sqrt{\lambda}+K\lambda  
-\mu)_{+}\|_{L_{d}}.
\end{equation}
For $\mu\geq
K\lambda+\mu_{\theta}\lambda^{1/2}$
the factor of $N_{\mu}$ in \eqref{8.29.4}  is dominated by
$$
N'\lambda^{-1/2}\|(|b|\sqrt{\lambda}  
-\mu_{\theta}\sqrt{\lambda})_{+}\|_{L_{d}}
=N' \|(|b| 
-\mu_{\theta})_{+}\|_{L_{d}}\leq N'\theta.
$$
  Thus, if $\theta$ is chosen in such a way that
$N'\theta\leq1/2$, then $N_{\mu}\leq 2N'\lambda^{-1}$,
which is equivalent to \eqref{9.10.6}
for $p=d$ as it was explained above. 

This proves the theorem for $p=d$.
For general $p\geq d$ it suffices to
 use the  Marcinkiewicz interpolation theorem
as in the proof of Theorem \ref{theorem 8.25.1}.
The theorem
is proved.

\end{document}